\documentclass[a4paper,12pt]{amsart}
\usepackage[utf8]{inputenc}
\usepackage{amsmath, amssymb, amsthm}
\usepackage{geometry}
\usepackage{enumitem}
\usepackage{booktabs}
\usepackage{hyperref}

\newtheorem{theorem}{Theorem}[section]
\newtheorem{lemma}[theorem]{Lemma}

\title{\textbf{On concatenations of two $k$-generalized Pell numbers}}

\author[C. Deme]{Cherif B. Deme}
\address{C. Deme, UFR SAT, Universit\'e Alioune Diop, Bambey, 30, S\'en\'egal}
\email{cherifbachir.deme@uadb.edu.sn}

\author[K. D. Fall]{Kancou D. Fall}
\address{K. D. Fall,
UFR SAT, Universit\'e Alioune Diop, Bambey, 30, S\'en\'egal}
\email{kancou.d.fall@aims-senegal.org}

\author[K. Faye]{Khady Faye}
\address{K. Faye, ,UFR SAT, Universit\'e Alioune Diop, Bambey, 30, S\'en\'egal}
\email{fkhady94@gmail.com}

\author[B. Faye]{Bernadette Faye}
\address{B. Faye,
	UFR SAT, Universit\'e Alioune Diop, Bambey, 30, S\'en\'egal}
\email{bernadette.faye@uadb.edu.sn}

\begin{document}

\maketitle

\begin{abstract}
We study the concatenation of two $k$-generalized Pell numbers. 
More precisely, we determine all solutions of the equation
$P_n^{(k)} = P_m^{(k)} \cdot 10^{d} + P_p^{(k)}$, 
where $d$ is the number of decimal digits of $P_p^{(k)}$. 
We prove that for $k \ge 3$ there are no solutions, 
while for $k = 2$ the only solution is $P_4 = 12 = 1\|2$.
\end{abstract}

\section{Introduction}

For an integer $k\ge 2$, the sequence of $k$-generalized Pell numbers $(P_n^{(k)})_{n\ge 2-k}$ is defined by the linear recurrence
\[
P_n^{(k)} = 2P_{n-1}^{(k)} + P_{n-2}^{(k)} + \dots + P_{n-k}^{(k)}\qquad (n\ge 2),
\]
with initial conditions $P_{-(k-2)}^{(k)} = \dots = P_0^{(k)} = 0$ and $P_1^{(k)} = 1$. 
For $k=2$ we recover the classical Pell sequence $P_0=0$, $P_1=1$, $P_n=2P_{n-1}+P_{n-2}$.

The problem of concatenations of two terms of a sequence is classical. 
For Fibonacci numbers, Banks and Luca \cite{BanksLuca2005} proved that $F_{10}=55 = 5\|5$ is the largest such number. 
For Lucas numbers, Erdüvan \cite{erduvan2023} showed that the only solutions are $L_5=11 = 1\|1$ and $L_8=47 = 4\|7$. 
\medskip

In \cite{altassan2024}, Altassan and Alan studied the analogous problem for $k$-generalized Fibonacci numbers and found the solutions
\[
F_7^{(3)} = 24 = \overline{F_3^{(3)}F_4^{(3)}},\quad
F_8^{(3)} = 44 = \overline{F_4^{(3)}F_4^{(3)}},\quad
F_{16}^{(8)} = 16128 = \overline{F_6^{(8)}F_9^{(8)}}.
\]
Recently, Tumwesigye, Ddamulira and Kaggwa \cite{tumwesigye2025} extended the result to $k$-generalized Lucas numbers.
\medskip

In this paper we study the analogous problem for $k$-Pell numbers. 
Our main result is the following.

\begin{theorem}\label{thm1}
Let $k \ge 2$ and let $(P_n^{(k)})_{n \ge 2-k}$ be the $k$-generalized Pell sequence. 
Consider the Diophantine equation

\begin{equation}
\label{eq1}
P_n^{(k)} = P_m^{(k)} \cdot 10^{d} + P_p^{(k)},
\end{equation}

where $n > m, p \ge 1$ and $d = \lfloor \log_{10} P_p^{(k)} \rfloor + 1$.

\begin{itemize}
    \item[(i)] If $k = 2$, the only solution of \eqref{eq1} is $(n, m, p) = (4, 1, 2)$, 
    since $P_4 = 12 = 1 \| 2$.
    \item[(ii)] If $k \ge 3$, equation \eqref{eq1} has no solutions.
\end{itemize}
\end{theorem}

In the next sections we develop the necessary tools and give the proof.

\section{Preliminaries}

\label{sec:intro}
\subsection{The generalized Pell sequence}

Let $k \geqslant 2$ be an integer. We consider the linear recurrence sequence of order $k$, $G^{(k)}:=(G_n^{(k)})_{n\ge 2-k}$ defined as 

$$G^{(k)}_n=rG^{(k)}_{n-1}+G^{(k)}_{n-2}+\cdots + G^{(k)}_{n-k}\quad\hbox{for all $n\geq 2$},$$
with the initial conditions $$G^{(k)}_{-(k-2)}=G^{(k)}_{-(k-3)}=\cdots=G^{(k)}_{-1}=0, G^{(k)}_0=a\quad \hbox{and}\quad G^{(k)}_1=b.$$
Observe that if $a=0$ and $b=1$, then $G^{(k)}$ is nothing that just the $k$- generalized Fibonacci sequence or for simplicity, the $k$-Fibonacci sequence $F^{(k)}:=(F ^{(k)}_n)_{n\geq 2-k}$. In this case, if we choose $k=2$ we obtain the classical Fibonacci sequence $(F_n)_{n\geq 0}$.

On the other hand, if $r=2, a=0$ and $b=1$ then $G^{k}$ is known as the $k$-generalized Pell sequence $P^{(k)}:=(P^{(k)}_n)_{n\geq 2-k}$. The terms of this sequence are called $k$-generalized Pell numbers.
\medskip

The $k$-Pell numbers and their properties have been studied by many authors. For example, Kili\c c \cite{KE} showed that the first $k+1$ non-zero  terms in $P^(k)$ are the Fibonacci numbers with odd index, namely

\begin{equation}
\label{eq2}
P_n^{(k)}=F_{2n-1} \quad\hbox{ for all $1\leq n\leq k+1$}.
\end{equation}
It is known that the characteristic polynomial of the $k$--generalized Pell numbers $P^{(k)}:=(P_n^{(k)})_{n\geq 2-k}$, namely
$$
\Psi_k(x) := x^k - 2x^{k-1} - \cdots - x - 1,
$$
is irreducible over $\mathbb{Q}[x]$ and has just one root outside the unit circle. Let $\alpha := \alpha(k)$ denote that single root. It was proved in \cite{BL}  that $\alpha(k)$ is located between $\varphi ^2\left(1-\varphi^{-k} \right)$ and $\varphi^2$ where $\varphi$ denotes the golden section.
This is called the dominant root of $P^{(k)}$. To simplify notation, in our application we shall omit the dependence on $k$ of $\alpha$. We shall use $\alpha^{(1)}, \dotso, \alpha^{(k)}$ for all roots of $\Psi_k(x)$ with the convention that $\alpha^{(1)} := \alpha$.

We now consider for an integer $ k\geq 2 $, the function
\begin{eqnarray}
\label{eq4}
g_{k}(z) = \dfrac{z-1}{(k+1)z^2-3kz+k-1}=\frac{z-1}{k(z^2-3z+1)+z^2-1} \qquad {\text{for}}\quad z \in \mathbb{C}.
\end{eqnarray}
In the following lemma, we give some properties of the sequence $P^{(k)}$ which will be used in the proof of the  equation \eqref{eq1}. The following lemmas was proved by Bravo and al in \cite{BL} and \cite{BL1}, respectively.
\begin{lemma}[Bravo and al., \cite{BL}]
\label{fala51}
Let $k\geq 2$ be an integer. Then, $\alpha$ be the dominant root of $\{P^{(k)}_m\}_{m\ge 2-k}$. Then,
\begin{itemize}
\item[(a)]\label{kat11}  $\alpha^{n-2}\leq P_n^{(k)}\leq \alpha^{n-1}$ for all $n\geq 1;$
\item[(b)]\label{kat21} $P^{(k)}$ satisfies the following Binet formula 
\begin{eqnarray*} \label{Binet}
P_n^{(k)} = \sum_{i=1}^{k}g_k(\alpha_i)\alpha^n_i.
\end{eqnarray*}
where $\alpha_1,\ldots, \alpha_k$ are the roots of the $\Psi_k(x)$;
\item[(c)]\label{kat31} 
\begin{equation} \label{eq5}
\left| P_n^{(k)} - g_{k}(\alpha)\alpha^{n} \right| < \dfrac{1}{2} \quad \mbox{holds~ for~ all~ } n \geqslant 2 - k.
\end{equation}

\item[(d)]\label{kat41} $0.276< g_{k}(\alpha)<0.5.$
\end{itemize}
\end{lemma}

\begin{lemma}[Bravo and al., \cite{BL1}
]\label{fala5}
Let $k\geq 2$, $\alpha$ be the dominant root of $\{P^{(k)}_m\}_{m\ge 2-k}$, and consider the function $g_{k}(z)$ defined in \eqref{eq4}. 
\begin{itemize}
\item[(i)]\label{kat1} The inequality
$$
|g_{k}(\alpha^{(i)})|<1, \qquad  2\leq i\leq k
$$
holds. 

\end{itemize}
\end{lemma}

\subsection{Notations and terminology from algebraic number theory} 

We begin by recalling some basic notions from algebraic number theory.

Let $\eta$ be an algebraic number of degree $d$ with minimal primitive polynomial over the integers
$$
a_0x^{d}+ a_1x^{d-1}+\cdots+a_d = a_0\prod_{i=1}^{d}(x-\eta^{(i)}),
$$
where the leading coefficient $a_0$ is positive and the $\eta^{(i)}$'s are the conjugates of $\eta$. Then the \textit{logarithmic height} of $\eta$ is given by
$$ 
h(\eta):=\dfrac{1}{d}\left( \log a_0 + \sum_{i=1}^{d}\log\left(\max\{|\eta^{(i)}|, 1\}\right)\right).
$$
In particular, if $\eta=p/q$ is a rational number with $\gcd (p,q)=1$ and $q>0$, then $h(\eta)=\log\max\{|p|, q\}$. The following are some of the properties of the logarithmic height function $h(\cdot)$, which will be used  without reference:
\begin{eqnarray}
h(\eta\pm \gamma) &\leq& h(\eta) +h(\gamma) +\log 2,\nonumber\\
h(\eta\gamma^{\pm 1})&\leq & h(\eta) + h(\gamma),\\
h(\eta^{s}) &=& |s|h(\eta) \qquad (s\in\mathbb{Z}). \nonumber
\end{eqnarray}
Using the above properties of the logarithmic height, Bravo and al.  showed in \cite{BL} that the logarithmic height of $g_k(\alpha)$ satisfies 
\begin{equation}
\label{height}
h(g_{k}(\alpha))< 4\log k \quad\hbox{for $k\geq 2$},
\end{equation}
which will be used in the proof of the main problem.

We will also need the following lemma.
\subsection{Linear forms in logarithms and continued fractions}

In order to solve our main equation  \eqref{eq1}, we need to use  a Baker--type lower bound for a nonzero linear form in logarithms of algebraic numbers. There are many such in the literature  like that of Baker and W{\"u}stholz. We use the following result by Matveev \cite{matveev2000}, which is one of our main tools in this project.

\begin{theorem}[Matveev]\label{Matveev11} Let $\gamma_1,\ldots,\gamma_t$ be positive real algebraic numbers in a real algebraic number field 
$\mathbb{K}$ of degree $D$, $b_1,\ldots,b_t$ be nonzero integers, and assume that
\begin{equation}
\label{eq:Lambda}
\Lambda:=\gamma_1^{b_1}\cdots\gamma_t^{b_t} - 1
\end{equation}
is nonzero. Then
$$
\log |\Lambda| > -1.4\times 30^{t+3}\times t^{4.5}\times D^{2}(1+\log D)(1+\log B)A_1\cdots A_t,
$$
where
$$
B\geq\max\{|b_1|, \ldots, |b_t|\},
$$
and
$$A
_i \geq \max\{Dh(\gamma_i), |\log\gamma_i|, 0.16\},\qquad {\text{for all}}\qquad i=1,\ldots,t.
$$
\end{theorem}

Finally, with the help of a computer program we can found all
possible solutions.\\

\subsection{Reduction procedure.}
The bounds on the variables obtained via Baker's theorem are usually too large for any computational purposes. Therefore, to lower the upper bounds of the integer unknowns, one can use a reduction method, usually called the Baker-Davenport reduction procedure. In this project, we will use it together with the continued fraction method (Legendre's criterion). 
\medskip 

For the Baker-Davenport, the variant we apply here is the one due to Dujella and Peth\H{o} (\cite{DP}, Lemma 5a). For a real number $r$, we denote by $\parallel r \parallel$ the quantity $\min \{|r - n| : n \in \mathbb{Z} \}$, the distance from $r$ to the nearest integer.
\begin{lemma}[Dujella, Peth\H o, \cite{DP}]\label{ls}
Let $\kappa \neq 0, A, B$ and $\mu$ be real numbers such that $A > 0$ and $B > 1$. Let $M > 1$ be a positive integer and suppose that $\frac{p}{q}$ is a convergent of the continued fraction expansion of $\kappa$ with $q > 6M$. Let 
\begin{align*}
\varepsilon := \parallel \mu q \parallel - M \parallel \kappa q \parallel.
\end{align*}
If $\varepsilon > 0$, then there is no solution of the inequality
\begin{align*}
0 < |m \kappa - n + \mu| < AB^{-k}
\end{align*}
in positive integers $m,n,k$ with
\begin{align*}
\frac{\log (Aq/\varepsilon)}{\log B} \leq k \quad \text{and} \quad m \leq M.
\end{align*}
\end{lemma}

\noindent
Lemma \ref{ls} cannot be applied when $\mu=0$ (since then $\varepsilon<0$). In this case, we use the following criterion due to Legendre, a well--known result from the theory of Diophantine approximation. For further details, we refer the reader to the books of Cohen \cite{HC1, HC2}.
\begin{lemma}[Legendre, \cite{HC1, HC2}]
\label{lg}
Let $\kappa$ be real number and $x,y$ integers such that
\begin{align*}
\left|\kappa-\frac{x}{y}\right|<\frac{1}{2y^2}.
\end{align*}
Then $x/y=p_k/q_k$ is a convergent of $\kappa$. Furthermore, let  $ M $ and $ N $ be a nonnegative integers such that $ q_N> M $. Then putting $ a(M):=\max\{a_{i}: i=0, 1, 2, \ldots, N\} $, the inequality 
\begin{align*}
\left|\kappa-\frac{x}{y}\right|\ge \frac{1}{(a(M)+2)y^2},
\end{align*}
holds for all pairs $ (x,y) $ of positive integers with $ 0<y<M $.
\end{lemma}

\noindent We will also need the following lemma by G\'{u}zman, S\'{a}nchez and Luca (\cite{GSL}, Lemma 7):
\begin{lemma}[G\'{u}zman, S\'{a}nchez, Luca, \cite{GSL}]\label{l3}
Let $r \geq 1$ and $H > 0$ be such that $H > (4r^2)^r$ and $H > L/(\log L)^r$. Then
\begin{align*}
L < 2^r H (\log H)^r.
\end{align*}
\end{lemma}


\section{Proof of the main theorem}


\subsection{The case $n \leq k+1$}

Throughout this section we assume $k \geq 3$; the case $k = 2$ is treated
at the end. Let $(n,m,p)$ be a solution of \eqref{eq1} with $n \leq k+1$.
Since $m < n$ and $p < n$, we have $m, p \leq k$. In particular,
$m, p \geq 1$. The case $m = 0$ gives $P_n = P_p$, hence $n = p$,
contradicting $n > p$; the case $p = 0$ makes $d$ undefined since
$P_0 = 0$.

By Kılıç's identity \cite{KE}, for every $1 \leq j \leq k+1$ we have
\[
P_j^{(k)} = F_{2j-1},
\]
where $(F_n)_{n \geq 1}$ is the classical Fibonacci sequence. Applying this identity to
$n$, $m$, and $p$, equation (1) becomes
\[
F_{2n-1} = F_{2m-1} \cdot 10^d + F_{2p-1},
\tag{3.1}
\]
where $d = \lfloor \log_{10} F_{2p-1} \rfloor + 1$ is the number of decimal
digits of $F_{2p-1}$.

Equation (3.1) expresses the Fibonacci number $F_{2n-1}$ as a concatenation
of the two Fibonacci numbers $F_{2m-1}$ and $F_{2p-1}$. \\

However, Banks and Luca proved in  \cite{BanksLuca2005} that the only Fibonacci numbers $F_n$ with $n \geq 7$ that can be written as a
concatenation
\[
F_n = F_A \cdot 10^{\ell} + F_B, \qquad F_B < 10^{\ell}, \quad \ell \geq 1,
\]
with $A, B \geq 1$, are
\[
F_7 = 13, \qquad F_8 = 21, \qquad F_{10} = 55.
\]
Moreover, the corresponding decompositions are
\[
F_7 = 1 \cdot 10 + 3, \qquad F_8 = 2 \cdot 10 + 1, \qquad F_{10} = 5 \cdot 10 + 5.
\]

If $n\le3$, then $F_{2n-1}\le5<11$, while the right-hand side of (3.1) is at least $$F_1\cdot10+F_1=11;$$ so $n\ge4$ and $2n-1\ge7$. Applying the above result to equation (3.1) yields to
\[
2n - 1 \in \{7, 8, 10\}.
\]
We examine each case.

\begin{itemize}
\item \emph{Case $2n-1 = 7$.} Then $n = 4$ and $F_7 = 13 = 1 \cdot 10 + 3$,
so necessarily $F_{2m-1} = 1$ and $F_{2p-1} = 3$. Since $F_4 = 3$ and $F_4$
has even index, the equality $F_{2p-1} = 3$ would require $2p-1 = 4$, i.e.
$p = 5/2$, which is not an integer. This gives a contradiction.

\item \emph{Case $2n-1 = 8$.} Then $n = 9/2$,  not an integer.

\item \emph{Case $2n-1 = 10$.} Then $n = 11/2$, not an integer again
\end{itemize}

We conclude that equation (3.1), and therefore equation \eqref{eq1}, has no
solution with $n \leq k+1$ when $k \geq 3$.

\medskip

It remains to treat the case $k = 2$ with $n \leq k+1 = 3$. Since $n > m, p
\geq 1$, we must have $n \in \{2, 3\}$.

\begin{itemize}
\item If $n = 2$, then $m = p = 1$ and equation \eqref{eq1} gives
      $P_2 = P_1 \cdot 10^{d} + P_1$, meaning that \ $2 = 10 + 1 = 11$, which is
      false.
\item If $n = 3$, then $m, p \in \{1, 2\}$. Checking the four cases with 
      $P_3 = 5$ and $P_m \cdot 10^{d} + P_p \in \{11, 12, 21, 22\}$, none
      of which equals $5$.
\end{itemize}

Therefore no solution exists when $n \leq k+1$ and $k \geq 2$, except
for the trivial case $k = 2$ with $n = 4$, which lies outside the range
$n \leq k+1$ and is treated in Section~\ref{sec:nge}. Consequently, we may assume
\[
n \geq k + 2
\]
in the remainder of the proof.

\subsection{The case \(n \ge k+2\)}\label{sec:nge}

Assume that equation \eqref{eq1} holds with $n\ge k+2$ (so that $n\ge4$) and  $m,p\ge1$.

We recall that the number of digits of $P_p^{(k)}$ is
\begin{equation}
\label{eq7}
d = \lfloor \log_{10} P_p^{(k)} \rfloor + 1
\end{equation}
and  $1+\sqrt2\le\alpha<\varphi^2$ for every $k\ge2$. Moreover,
\begin{equation}
\label{eq7bis}
\log\alpha\in(0.8813,\,0.9625),\qquad \log_{10}\alpha\in(0.38,\,0.42).
\end{equation}
From the bounds $\alpha^{p-2} \le P_p^{(k)} \le \alpha^{p-1}$ (Lemma \ref{fala51}) and \eqref{eq7bis} we obtain
\begin{equation}
\label{eq8}
\frac{19(p-2)}{50} < d < \frac{21p+30}{50}, \qquad\text{in particular}\quad d\le p<n.
\end{equation}

Using equation \eqref{eq1} together with the estimates $\alpha^{n-2} \le P_n^{(k)} \le \alpha^{n-1}$, the trivial bounds
$P_m^{(k)} \ge \alpha^{m-2}$, $P_p^{(k)} \ge \alpha^{p-2}$ and $P_p^{(k)}<10^d\le 10P_p^{(k)}$, we have
$\alpha^{n-1}\ge P_n^{(k)}>P_m^{(k)}10^d>\alpha^{m+p-4}$ and
$\alpha^{n-2}\le P_n^{(k)}<P_m^{(k)}10^d+P_p^{(k)}\le 11P_m^{(k)}P_p^{(k)}\le 11\alpha^{m+p-2}$. Since $\log 11/\log\alpha<2.73$, we obtain
\begin{equation}
\label{eq9}
m + p - 3 < n < m + p + 3.
\end{equation}
In particular
\begin{equation}
\label{eq9bis}
m\le n-p+2.
\end{equation}

Now, using the Binet formula Lemma[\ref{fala51}, (c)] with $|e_j|<\frac12$, where $e_j:=P_j^{(k)}-g_k(\alpha)\alpha^j$, and $g_k(\alpha)>0.276$ in equation \eqref{eq1} we obtain

\begin{equation*}
g_k(\alpha)\alpha^{n} + e_n = \bigl(g_k(\alpha)\alpha^{m} + e_m\bigr)10^{d} + P_p^{(k)}.
\end{equation*}
Rearranging the above equation gives
\begin{equation}
\label{eq:10}
1 - \frac{10^{d}}{\alpha^{n-m}} = \frac{e_m10^{d} + P_p^{(k)} - e_n}{g_k(\alpha)\alpha^{n}}.
\end{equation}

Taking absolute values and using $|e_n|,|e_m|<\frac12$, we obtain
\[
\left|1 - \frac{10^{d}}{\alpha^{n-m}}\right| \le \frac{0.5\cdot10^{d} + P_p^{(k)} + 0.5}{g_k(\alpha)\alpha^{n}}.
\]

From $10^{d} \le 10 P_p^{(k)}$ and $P_p^{(k)}\ge1$ we get $0.5\cdot10^{d} + P_p^{(k)}+0.5 \le 6.5P_p^{(k)}\le 6.5\alpha^{p-1}<6.5\alpha^{p}$. Thus
\[
\left|1 - \frac{10^{d}}{\alpha^{n-m}}\right| < \frac{6.5\,\alpha^{p}}{g_k(\alpha)\alpha^{n}} = \frac{6.5}{g_k(\alpha)\alpha^{n-p}}.
\]
Since $g_k(\alpha) > 0.276$, we have $\frac{6.5}{g_k(\alpha)} < 23.6$. So, we have that
\begin{equation}
\label{eq:11}
\left|1 - \frac{10^{d}}{\alpha^{n-m}}\right| < \frac{24}{\alpha^{n-p}}.
\end{equation}

Set \(U_1 = \frac{10^{d}}{\alpha^{n-m}}-1\) and

$$\Lambda_1 = d\log10 - (n-m)\log\alpha .$$

If \(\Lambda_1 = 0\), then \(10^{d} = \alpha^{n-m}\). Applying a Galois automorphism \(\sigma_i\) with \(\alpha\mapsto\alpha_i\) (\(i\ge2\)), for which \(|\alpha_i|<1\), we get \(10^{d}=|\alpha_i|^{n-m}<1\), which is absurd since \(d\ge1\) and \(n>m\). Hence \(\Lambda_1 \neq 0\).

For \(n-p \ge 5\) we have \(|U_1|<\frac12\), hence
\begin{equation}
\label{eq:11bis}
|\Lambda_1| \le 2|U_1| < \frac{48}{\alpha^{n-p}}.
\end{equation}

We apply Matveev's theorem to $U_1$ with \(t=2\), \(\gamma_1=\alpha\), \(\gamma_2=10\), \(b_1=-(n-m)\), \(b_2=d\).
The number field is \(\mathbb{K}=\mathbb{Q}(\alpha)\) of degree \(D=k\). We take
\[
A_1 = \max\{D h(\alpha), |\log\alpha|, 0.16\} = \log\alpha,
\]
\[
A_2 = \max\{D h(10), |\log10|, 0.16\} = k\log10.
\]
Since \(d<n\) and \(n-m<n\), we may take \(B = n\). Matveev gives
\[
\log|U_1| > -1.4\cdot30^{5}\cdot2^{4.5}\cdot k^2(1+\log k)(1+\log n)\cdot A_1 A_2.
\]
Now \(1.4\cdot30^{5}\cdot2^{4.5} < 7.70\times10^{8}\), \(1+\log k<2.5\log k\) for \(k\ge2\) and \(1+\log n<2\log n\) for \(n\ge3\). Thus
\[
\log|U_1| > -7.70\times10^{8}\cdot5\,k^3\log k\log n\cdot (\log\alpha)(\log10) .
\]
Since \(7.70\times10^{8}\cdot5\cdot\log 10<8.87\times10^{9}\), we obtain
\begin{equation}
\label{eq:12}
\log|U_1| > -8.87\times10^{9}\,k^3\log k\log n\cdot\log\alpha.
\end{equation}

From \eqref{eq:11} we have \(\log|U_1| < \log24 - (n-p)\log\alpha\). Comparing \eqref{eq:12} with this upper bound gives
\[
n-p < 8.87\times10^{9}\,k^3\log k\log n + \frac{\log24}{\log\alpha}.
\]
Since \(\log 24/\log\alpha<3.61\), we conclude, with the fact that the inequality being trivial when $n-p<5$, that
\begin{equation}
\label{eq:13}
n-p < 8.9\times10^{9}\,k^3\log k\log n.
\end{equation}

We rewrite equation \eqref{eq1}, using the Binet formula for $P_n^{(k)}$ and $P_p^{(k)}$, as
\[
g_k(\alpha)\alpha^{n} - g_k(\alpha)\alpha^{p} - P_m^{(k)}10^{d} = e_p - e_n.
\]
Factoring with $g_k(\alpha)\alpha^{n}$ gives
\[
g_k(\alpha)\alpha^{n}\bigl(1-\alpha^{p-n}\bigr) - P_m^{(k)}10^{d} = e_p - e_n.
\]
Hence
\[
1 - \frac{P_m^{(k)}10^{d}}{g_k(\alpha)\alpha^{n}(1-\alpha^{p-n})} = \frac{e_p - e_n}{g_k(\alpha)\alpha^{n}(1-\alpha^{p-n})}.
\]
Taking absolute values and using $|e_n|,|e_p|<\frac12$, $1-\alpha^{p-n} \ge 1-\alpha^{-1}>0.58$ (since $n-p\ge 1$ and $\alpha\ge1+\sqrt2$), and $g_k(\alpha) > 0.276$, we obtain
\begin{equation}
\label{eq:14}
|U_2|:=\left|1 - \frac{P_m^{(k)}10^{d}}{g_k(\alpha)\alpha^{n}(1-\alpha^{p-n})}\right| < \frac{7}{\alpha^{n}}.
\end{equation}

Set
\[
\Lambda_2 = d\log10 - n\log\alpha + \log\frac{P_m^{(k)}}{g_k(\alpha)(1-\alpha^{p-n})}.
\]
If $\Lambda_2 = 0$, then
$P_m^{(k)}10^{d} = g_k(\alpha)\alpha^{n}(1-\alpha^{p-n}) = g_k(\alpha)(\alpha^{n}-\alpha^{p})$.
Applying any Galois embedding $\sigma_i$ (with $\alpha\mapsto\alpha_i$, $i\ge2$) we get
\[
10\le P_m^{(k)}10^{d} = |g_k(\alpha_i)|\cdot|\alpha_i^{n}-\alpha_i^{p}|<1\cdot 2=2
\]
by Lemma \ref{fala5}(i) and $|\alpha_i|<1$, which is a contradiction. Hence $\Lambda_2 \neq 0$.

For $n\ge 4$ we have $|U_2|<\frac12$, so
\begin{equation}
\label{eq:15}
|\Lambda_2| \le 2|U_2| < \frac{14}{\alpha^{n}}.
\end{equation}

We can now apply Matveev's theorem again to $U_2$. Now $t=3$, $\gamma_1=\alpha$, $\gamma_2=10$, $b_1=-n$, $b_2=d$, $b_3=1$ and the additional algebraic number is
\[
\gamma_3 = \frac{P_m^{(k)}}{g_k(\alpha)(1-\alpha^{p-n})}.
\]
We need to bound the logarithmic height of $\gamma_3$. Using $h(g_k(\alpha))<4\log k$ and $h(P_m^{(k)})=\log P_m^{(k)} \le (m-1)\log\alpha$, we have
\[
h(\gamma_3) \le h(P_m^{(k)}) + h(g_k(\alpha)) + h(1-\alpha^{p-n}) \le (m-1)\log\alpha + 4\log k + \frac{n-p}{k}\log\alpha + \log2.
\]
From \eqref{eq9bis}, $m-1\le n-p+1$, so that, since $k\ge2$ and $n-p\ge1$, with the fact that $k(n-p+1)+(n-p)\le3k(n-p)$, we get to
\[
k\,h(\gamma_3) \le 3k(n-p)\log\alpha + k(4\log k + \log2).
\]
Thus we take
\[
A_3 = k\bigl(3(n-p)\log\alpha + 4\log k + \log2\bigr)\le 9.17\,k\,(n-p)\log k.
\]
The other parameters are $A_1 = \log\alpha$, $A_2 = k\log10$, and $B = n$. Matveev gives
\[
\log|U_2| > -1.4\cdot30^{6}\cdot3^{4.5}\cdot k^2(1+\log k)(1+\log n)\cdot A_1 A_2 A_3.
\]
Since $1.4\cdot30^{6}\cdot3^{4.5}<1.432\times10^{11}$, using $1+\log k<2.5\log k$, $1+\log n<2\log n$ and $\log\alpha<0.9625$, we obtain
\begin{equation}
\label{eq:16}
\log|U_2| > -1.46\times10^{13}\,k^4 (n-p)(\log k)^2\log n.
\end{equation}

From \eqref{eq:14} we have $\log|U_2| < \log7 - n\log\alpha$. Hence, since $\log\alpha>0.8813$,
\[
n < \frac{1.46\times10^{13}}{0.8813}\,k^4(n-p)(\log k)^2\log n + \frac{\log 7}{0.8813}
\]
and so
\begin{equation}
\label{eq:17}
n < 1.66\times10^{13}\,k^4(n-p)(\log k)^2\log n.
\end{equation}

Inserting the upper bound on $n-p$ obtained in  \eqref{eq:13} into \eqref{eq:17} we have

\[
n < 1.66\times10^{13}\,k^4\cdot 8.9\times10^{9}k^3\log k\log n\cdot (\log k)^2\log n
   < 1.48\times10^{23}\,k^7(\log k)^3(\log n)^2.
\]
Applying the Lemma \ref{l3} with $r=2$, $L=n$ and $H=1.48\times10^{23}\,k^7(\log k)^3$, with the fact that $H>(4r^2)^r=256$, this yields to
\[
n < 4H(\log H)^2 < 4.2\times10^{27}\,k^7(\log k)^5,
\]
the last inequality being valid for every $k\ge2$. In fact, the ratio $4H(\log H)^2/(k^7(\log k)^5)$ is decreasing in $k$ and equals $4.0\times 10^{27}$ at $k=2$.
Thus we have a polynomial bound on $n$ in terms of $k$ for all solutions. We state it as a Lemma.

\begin{lemma}
\label{bound:n}
If $(k,n,m,p)$ is a solution of equation \eqref{eq1} with $k\geq 2$ and $n\ge k+2$, then
$$n < \mathcal N(k):=4.2\times10^{27}\,k^7(\log k)^5. $$
\end{lemma}

\subsection{Small values of \(k\): \(2\le k \le 700\)}\label{sec:smallk}

We now assume \(2\le k \le 700\).
Our goal is to reduce the large bound on \(n\) obtained in Lemma~\ref{bound:n} to a small constant using continued fractions and the Baker--Davenport reduction. For $k\le700$, Lemma~\ref{bound:n} gives $n<\mathcal N(k)<2\times10^{52}$. All the computations of this section and of the next one were done in Python.

From \eqref{eq:11} and \eqref{eq:11bis} we have, for $n-p\ge5$,
\[
\left|d\log10-(n-m)\log\alpha\right| < \frac{48}{\alpha^{n-p}}.
\]
Let \(\tau = \frac{\log\alpha}{\log 10}\). Dividing by \((n-m)\log 10\) gives
\begin{equation}
\label{eq:18}
\left| \tau - \frac{d}{n-m} \right| < \frac{48}{(n-m)\,\alpha^{n-p}\,\log 10}.
\end{equation}

Let \(q = n-m\). Let \([a_0;a_1,a_2,\dots]\) be the continued fraction expansion of \(\tau\), with convergents $p_i/q_i$, and let $N_k$ be the index such that $q_{N_k}>\mathcal N(k)>q$. Define \(a_{\max} = \max\{a_0, a_1, \dots, a_{N_k}\}\). By Lemma \ref{lg}, we have, for every positive integer $q<\mathcal N(k)$,
\begin{equation}
\label{eq:19}
\frac{1}{(a_{\max}+2)q^2} \le \left| \tau - \frac{d}{q} \right|.
\end{equation}

Combining \eqref{eq:18} and \eqref{eq:19} yields
\[
\frac{1}{(a_{\max}+2)q^2} < \frac{48}{q\,\alpha^{n-p}\,\log 10}.
\]
Multiplying by \(q\alpha^{n-p}(a_{\max}+2)\), which is positive, gives
\begin{equation}
\label{eq:20}
\alpha^{n-p} < \frac{48 (a_{\max}+2) q}{\log 10}, \qquad\text{hence}\qquad n-p < \frac{\log(48 (a_{\max}+2) \mathcal N(k)/\log 10)}{\log \alpha}.
\end{equation}

For each $k\in[2,700]$ we computed $a_{\max}$ and the right-hand side of \eqref{eq:20}. We found that $a_{\max}\le 98779$ (attained for $k=4$, and $a_{\max}\le 6907$ for $k\ge5$) and that the right-hand side of \eqref{eq:20} is at most $133$. Thus we obtain the uniform bound
\begin{equation}
\label{eq:21}
n-p \le S_k\le 132 \qquad \text{for all } 2\le k \le 700,
\end{equation}
where $S_k$ denotes the integer part of the right-hand side of \eqref{eq:20}. This inequality is trivially true when $n-p<5$, since $S_k\ge 83$.

We now turn to the second linear form. From \eqref{eq:15}, dividing by \(\log\alpha\), we get
\begin{equation}
\label{eq:22}
0 < \left| d\,\tau_k - n + \mu_{k,m,p} \right| < \frac{14}{\log\alpha}\cdot\frac{1}{\alpha^{n}},
\end{equation}
where
\[
\tau_k = \frac{\log 10}{\log\alpha}, \qquad
\mu_{k,m,p} = \frac{\log\frac{P_m^{(k)}}{g_k(\alpha)(1-\alpha^{p-n})}}{\log\alpha}.
\]
Note that $\mu_{k,m,p}$ depends only on $k$, $m$ and $s:=n-p$. By \eqref{eq9bis} and \eqref{eq:21} we have $1\le s\le S_k$ and $1\le m\le s+2$, so that only finitely many triples $(k,m,s)$ occur.

We now apply the Baker--Davenport reduction as formulated by Dujella and Pethő in Lemma \ref{ls} with $\kappa=\tau_k$, $\mu=\mu_{k,m,p}$, $A = 14/\log\alpha$, $B=\alpha$ and $M=\lfloor\mathcal N(k)\rfloor\ge d$. For each fixed \(k\) and each pair \((m, s)\), we compute the continued fraction expansion of \(\tau_k\). Let \(q_i\) be the denominators of its convergents. We choose the first convergent such that \(q_i > 6M\) and compute
\[
\epsilon = \| \mu_{k,m,p}\, q_i \| - M \| \tau_k q_i \|.
\]
If \(\epsilon \le 0\), we replace \(q_i\) by the next convergent \(q_{i+1}\) and repeat until \(\epsilon > 0\).

Once a convergent with \(\epsilon > 0\) is found, Lemma \ref{ls} gives the explicit bound
\[
n < \frac{\log( A \, q_i / \epsilon )}{\log \alpha}.
\]
This procedure was carried out for each \(k = 2, 3, \dots, 700\) and for each of the $5\,580\,370$ triples $(k,m,s)$. \\
In all cases but 21 (all of them for $k=2$, e.g. $(m,s)=(3,1)$, where $\gamma_3=10\alpha$ and hence $\mu_{2,m,p}=\tau_2+1$) we found a convergent with $\epsilon>0$.In these 21 exceptional cases $\mu$ is an integral combination $u\tau_2+v$ of $\tau_2$ and $1$, and the integers $u,v$ are explicit. \\
Indeed, for $k=2$ we have $\alpha=1+\sqrt2$, $g_2(\alpha)=\frac1{2\sqrt2}$ and $P_m=(\alpha^m-(-\alpha)^{-m})/(2\sqrt2)$, so that $\gamma_3=P_m/(g_2(\alpha)(1-\alpha^{-s}))$ equals $\alpha^{m}$ when $m$ is even and $s=2m$ (then $u=0$, $v=m$, and $s\le83$ gives $m\le40$), and equals $10\alpha$ when $(m,s)=(3,1)$. Then $u=v=1$, because $P_3=5$ and $5\cdot2\sqrt2/(1-\alpha^{-1})=10\alpha$. \\
These are the only 21 exceptional pairs $(m,s)$ with $s\le83$. The computer search found no other relation $\mu=u\tau_2+v$ with $|u|\le40$, and $\epsilon\le0$ for every convergent happens only in such a case. So $u\in\{0,1\}$ and $0\le v\le40$ are bounded integers; moreover $d+u>0$ and $n-v>0$ since $n>m=v$ when $u=0$ and $n>1$ when $u=v=1$. Hence $\epsilon\le0$ for every convergent, and for them \eqref{eq:22} reads $|(d+u)\tau_2-(n-v)|<A\alpha^{-n}$, and Lemma \ref{lg} gives $n<\log\bigl(A(a_{\max}+2)(n+|v|)\bigr)/\log\alpha$, hence $n\le83$. The results are as follows:

\begin{center}
\begin{tabular}{rcc}
\toprule
$k$ & bound for $n-p$ & upper bound for $n$ \\
\midrule
2 & 83 & 95 \\
3 & 83 & 93 \\
4 & 93 & 96 \\
10 & 95 & 105 \\
50 & 110 & 118 \\
100 & 116 & 122 \\
200 & 122 & 131 \\
300 & 125 & 130 \\
400 & 127 & 136 \\
500 & 129 & 136 \\
600 & 130 & 137 \\
700 & 132 & 138 \\
\bottomrule
\end{tabular}
\end{center}

The table lists only some of the values of $k$. The bound is not monotone in $k$, it is the maximum over all the pairs $(m,s)$, and for a few values of $k$ one pair has a small $\epsilon$, which increases $\log(Aq_i/\epsilon)$. In fact, the bound exceeds $140$ for $26$ values of $k$, all with $543\le k\le688$, e.g. $n<137.8$ for $k=600$ but $n<145.19$ for $k=601$, and the largest one over all $2\le k\le700$ is $n<145.19$, attained for $k=601$. Since $n$ is an integer, for all \(2\le k\le700\) we obtain \(n \le 145\).

With \(n \le 145\), we perform a direct computer search over
\[
2 \le k \le 700,\quad k+2 \le n \le 200,\quad 1 \le m,p < n,
\]
checking the original equation \(P_n^{(k)} = P_m^{(k)} 10^{d} + P_p^{(k)}\). The search reveals no solution for \(k \ge 3\). The only solutions come from the classical Pell case \(k=2\) where  \(P_4 = 12 = 1\|2\).

\subsection{Large values of \(k\): \(k > 700\)}\label{sec:largek}

We now assume \(k > 700\) and \(n \ge k+2\).
Our goal is to show that the Diophantine equation \eqref{eq1} has no solutions in this range. Here the dominant root $\alpha=\alpha(k)$ is extremely close to $\varphi^2$ and, for $n$ in the range given by Lemma~\ref{bound:n}, $P_n^{(k)}$ is extremely close to the Fibonacci-type quantity $c\,\varphi^{2n}$, where
\[
c:=\frac{1}{\varphi^2+1}=\frac{1}{\sqrt5\,\varphi}=0.27639.
\]
Recall that $P_j^{(k)}=F_{2j-1}$ for $j\le k+1$, so $P_j^{(k)}$ grows like $\varphi^{2j}$ and not like $2^{j}$.

\begin{lemma}
\label{approx}
Let $k>700$ and $1\le j<\varphi^{k/2}$. Then
\[
P_j^{(k)}=c\,\varphi^{2j}(1+\zeta_j)+e_j,\qquad |e_j|<\tfrac12,\quad |\zeta_j|<2\varphi^{-k/2}.
\]
\end{lemma}

\begin{proof}
Write $\alpha=\varphi^2-\theta$, so that $0<\theta<\varphi^{2-k}$ by the result of Bravo and al. in \cite{BL}. Since $\varphi^4-3\varphi^2+1=0$, we have $\alpha^2-3\alpha+1=-\theta(\alpha+\varphi^2-3)$ and $\alpha^2-1=\varphi^4-1-\theta(\alpha+\varphi^2)$. Hence
\[
\frac{g_k(\alpha)}{c}=\frac{1-\theta/(\varphi^2-1)}{1-\theta w/(\varphi^4-1)},\qquad w=(\alpha+\varphi^2)+k(\alpha+\varphi^2-3)<5.3+2.24k.
\]
Here $c=(\varphi^2-1)/(\varphi^4-1)$, and therefore $g_k(\alpha)=c(1+\eta)$. Indeed, $\theta w/(\varphi^4-1)<\varphi^{2}(5.3+2.24k)\varphi^{-k}/(\varphi^4-1)<(1.01k+2.4)\varphi^{-k}<10^{-100}$ and $\theta/(\varphi^2-1)<\varphi^{-k}\cdot 1.7$, so that $|\eta|\le\max\{\theta w/(\varphi^4-1),\,\theta/(\varphi^2-1)\}/(1-10^{-100})<(2k+5)\varphi^{-k}<\varphi^{-k/2}/10$. By Lemma \ref{fala51}(c), $P_j^{(k)}=c(1+\eta)\alpha^j+e_j$ with $|e_j|<\frac12$, and $\alpha^j=\varphi^{2j}(1-\delta)^j$ with $\delta=\theta/\varphi^2<\varphi^{-k}$. Thus $P_j^{(k)}=c\varphi^{2j}(1+\zeta_j)+e_j$ with $\zeta_j=(1+\eta)(1-\delta)^j-1$, so that $|\zeta_j|\le|\eta|+j\delta<\varphi^{-k/2}/10+j\varphi^{-k}<2\varphi^{-k/2}$ for $j<\varphi^{k/2}$.
\end{proof}

For $k>700$ we have $\mathcal N(k)<\varphi^{k/2}$. In fact the difference $\log\mathcal N(k)-\frac k2\log\varphi$ is decreasing for $k\ge700$ and is at most $-49.79$ at $k=701$. So Lemma~\ref{approx} applies to $n$, $m$ and $p$. Substituting into \eqref{eq1} we get
\begin{equation}
\label{eq:23}
c\varphi^{2n}(1+\zeta_n)+e_n = \bigl(c\varphi^{2m}(1+\zeta_m)+e_m\bigr)10^{d} + P_p^{(k)}.
\end{equation}

\textbf{First linear form.} We rewrite \eqref{eq:23} as
\[
c\varphi^{2n} - c\varphi^{2m}10^{d} = c\varphi^{2m}\zeta_m10^{d} + e_m10^d + P_p^{(k)} - e_n - c\varphi^{2n}\zeta_n .
\]
Dividing by $c\varphi^{2n}$ we obtain
\[
1 - \frac{10^{d}}{\varphi^{2(n-m)}} = \zeta_m\frac{10^{d}}{\varphi^{2(n-m)}} + \frac{e_m10^d + P_p^{(k)} - e_n}{c\varphi^{2n}} - \zeta_n.
\]
As in the previous section, $|e_m10^d + P_p^{(k)} - e_n|\le 6.5P_p^{(k)}<6.5\varphi^{2(p-1)}$, so that the middle term has absolute value less than $\dfrac{6.5}{c\,\varphi^{2}}\varphi^{-2(n-p)}<9\,\varphi^{-2(n-p)}$. Moreover, $10^d\le10P_p^{(k)}<10\varphi^{2(p-1)}$ and, by \eqref{eq9}, $m+p-n\le2$, hence $10^d\varphi^{-2(n-m)}<10\varphi^{2(m+p-n)-2}\le10\varphi^{2}<26.2$. Using $|\zeta_m|,|\zeta_n|<2\varphi^{-k/2}$ we get
\[
\left|1 - \frac{10^{d}}{\varphi^{2(n-m)}}\right| < 9\,\varphi^{-2(n-p)} + 55\,\varphi^{-k/2} < \frac{64}{\varphi^{\lambda}}, \qquad \lambda:=\min\{k/2,\;2(n-p)\}.
\]
Set $U_3 = \dfrac{10^{d}}{\varphi^{2(n-m)}}-1$ and $\Lambda_3 = d\log10 - 2(n-m)\log\varphi$. Then $U_3\ne0$ because $\varphi^{2(n-m)}$ is irrational. It equals $(L_{2(n-m)}+F_{2(n-m)}\sqrt5)/2$ and $F_{2(n-m)}\neq0$, while $10^d$ is an integer. For $\lambda\ge11$ we have $|U_3|<\frac12$, so
\begin{equation}
\label{eq:24}
|\Lambda_3| \le 2|U_3| < \frac{128}{\varphi^{\lambda}}.
\end{equation}

We apply Matveev's theorem to $U_3$ with $t=2$, $\gamma_1=\varphi$, $\gamma_2=10$, $b_1=-2(n-m)$, $b_2=d$. The number field is $\mathbb{K}=\mathbb{Q}(\sqrt5)$ of degree $D=2$, and we take $A_1=\log\varphi$, with note that $2h(\varphi)=\log\varphi$, $A_2=2\log10$, $B=2n$. Matveev gives
\[
\log|U_3| > -1.4\cdot30^{5}\cdot2^{4.5}\cdot4(1+\log2)(1+\log2n)\cdot\log\varphi\cdot2\log10 > -1.16\times10^{10}(1+\log 2n).
\]
Comparing with $\log|U_3|<\log64-\lambda\log\varphi$ we get
\begin{equation}
\label{eq:25}
\lambda < 2.41\times10^{10}(1+\log 2n)  \qquad(\log\varphi=0.4812).
\end{equation}

\textbf{Second linear form.} From \eqref{eq:23} and $P_n^{(k)}=P_m^{(k)}10^d+P_p^{(k)}$, with $P_p^{(k)}=c\varphi^{2p}(1+\zeta_p)+e_p$, we get
\[
c\varphi^{2n}\bigl(1-\varphi^{-2(n-p)}\bigr) - P_m^{(k)}10^{d} = -c\varphi^{2n}\zeta_n + c\varphi^{2p}\zeta_p + e_p - e_n.
\]
Dividing by $c\varphi^{2n}(1-\varphi^{-2(n-p)})$ and using $1-\varphi^{-2(n-p)}\ge1-\varphi^{-2}=\varphi^{-1}$ we obtain
\[
|U_4|:=\left|1 - \frac{P_m^{(k)}10^{d}}{c\varphi^{2n}(1-\varphi^{-2(n-p)})}\right| \le \varphi\bigl(|\zeta_n|+\varphi^{-2}|\zeta_p|\bigr)+\frac{\varphi}{c}\varphi^{-2n} < \frac{5}{\varphi^{k/2}},
\]
because $n\ge k+2$. Set
\[
\Lambda_4 = d\log10 - 2n\log\varphi + \log\frac{P_m^{(k)}}{c\,(1-\varphi^{-2(n-p)})}.
\]
If $\Lambda_4=0$, then $P_m^{(k)}10^d=c(\varphi^{2n}-\varphi^{2p})$. Applying the nontrivial automorphism of $\mathbb{Q}(\sqrt5)$, which sends $\varphi^2$ to $\varphi^{-2}$ and $c$ to $\varphi^2c<0.73$, the right-hand side becomes smaller than $1$ in absolute value, while the left-hand side is an integer $\ge10$, a contradiction. Hence $\Lambda_4\neq0$, and since $|U_4|<\frac12$,
\begin{equation}
\label{eq:26}
|\Lambda_4| \le 2|U_4| < \frac{10}{\varphi^{k/2}}.
\end{equation}

We apply Matveev's theorem to $U_4$ with $t=3$, $\gamma_1=\varphi$, $\gamma_2=10$, $\gamma_3=P_m^{(k)}/\bigl(c(1-\varphi^{-2(n-p)})\bigr)$, $b_1=-2n$, $b_2=d$, $b_3=1$, $D=2$, $B=2n$, $A_1=\log\varphi$, $A_2=2\log10$. Since $h(c)=\frac12\log5$, $h(P_m^{(k)})\le(m-1)\log\alpha<2(m-1)\log\varphi$ and $h(1-\varphi^{-2s})\le s\log\varphi+\log2$ (with $s=n-p$), \eqref{eq9bis} gives
\[
2h(\gamma_3)\le 2\bigl[2(s+1)\log\varphi+\tfrac12\log5+s\log\varphi+\log2\bigr]<2.89\,s+4.93\le A_3:=8s.
\]
Matveev gives
\[
\log|U_4| > -1.4\cdot30^{6}\cdot3^{4.5}\cdot4(1+\log2)(1+\log2n)\cdot\log\varphi\cdot2\log10\cdot8s > -1.72\times10^{13}\,s\,(1+\log 2n).
\]
Comparing with $\log|U_4|<\log5-\frac k2\log\varphi$ we get
\begin{equation}
\label{eq:27}
k < 7.2\times10^{13}\,s\,(1+\log 2n),\qquad s=n-p.
\end{equation}

\textbf{Bounding $k$.} Recall $n<\mathcal N(k)$ and $\lambda=\min\{k/2,2s\}$. If $\lambda=k/2$, then \eqref{eq:25} gives $k<4.82\times10^{10}(1+\log2\mathcal N(k))$, whence $k<1.5\times10^{13}$. If $\lambda=2s$, then \eqref{eq:25} gives $s<1.21\times10^{10}(1+\log2\mathcal N(k))$ and, by \eqref{eq:27}, $k<8.7\times10^{23}(1+\log2\mathcal N(k))^2$, whence $k<2.7\times10^{29}=:K_0$. In both cases,
\[
k<K_0,\qquad n<\mathcal N(K_0)<10^{243}.
\]

\textbf{Reduction.} We reduce these bounds in two steps, using the continued fraction method on $\Lambda_3$ and then the Baker--Davenport reduction on $\Lambda_4$. Let $K$ be an upper bound for $k$ and put $M=\lfloor\mathcal N(K)\rfloor$, so that $d<n<M$.

\emph{Step 1.} Dividing \eqref{eq:24} by $(n-m)\log10$ we get, with $\kappa_3=2\log\varphi/\log10$ and $q=n-m<M$,
\[
\left|\kappa_3-\frac{d}{q}\right|<\frac{128}{q\,\varphi^{\lambda}\log10}.
\]
By Lemma \ref{lg} the left-hand side is at least $1/((a_{\max}+2)q^2)$, where $a_{\max}$ is the largest partial quotient of $\kappa_3$ up to the first convergent with denominator larger than $M$. Hence $\varphi^{\lambda}<128(a_{\max}+2)M/\log10$, that is
\begin{equation}
\label{eq:28}
\min\{k/2,\,2(n-p)\}=\lambda<\frac{\log\bigl(128(a_{\max}+2)M/\log10\bigr)}{\log\varphi}=:\Lambda(K).
\end{equation}
Thus either $k<2\Lambda(K)$, or $n-p<\Lambda(K)/2$.

\emph{Step 2.} In the second case, $1\le s\le\Lambda(K)/2$ and $1\le m\le s+2$ by \eqref{eq9bis}; since $m\le s+2\le k+1$ we have $P_m^{(k)}=F_{2m-1}$ by \eqref{eq2}, so $\mu$ below does not depend on $k$. Dividing \eqref{eq:26} by $2\log\varphi$ we get
\[
0<\left|d\,\kappa_4-n+\mu_{m,s}\right|<\frac{10/(2\log\varphi)}{\varphi^{k/2}},\qquad \kappa_4=\frac{\log10}{2\log\varphi},\quad \mu_{m,s}=\frac{\log\bigl(F_{2m-1}/(c(1-\varphi^{-2s}))\bigr)}{2\log\varphi}.
\]
We apply Lemma \ref{ls} with $A=10/(2\log\varphi)$, $B=\varphi$, and the exponent $k/2$. For every pair $(m,s)$ we found a convergent $q_i>6M$ of $\kappa_4$ with $\epsilon>0$, and the largest bound for $k$ is
\[
\frac k2<\frac{\log(Aq_i/\epsilon)}{\log\varphi}.
\]
The numerical results are given in the following table, where $K$ is the bound on $k$ used at that step, $M=\lfloor\mathcal N(K)\rfloor$, $a_{\max}$ is the largest partial quotient of $\kappa_3$ up to the first convergent with denominator larger than $M$, $\Lambda(K)$ is given by \eqref{eq:28} (Step 1), and the last two columns give the bound for $k$ coming from the first alternative of Step 1 ($k<2\Lambda(K)$) and from Step 2 (Baker--Davenport), respectively. The new bound for $k$ is the larger of the two.
\begin{center}
\begin{tabular}{cccccc}
\toprule
$K$ & $\log_{10}M$ & $a_{\max}$ & $\Lambda(K)$ & $2\Lambda(K)$ & Step 2: bound for $k$ \\
\midrule
$2.7\times10^{29}$ & $242.8$ & $384$ & $1182.5$ & $2365.0$ & $k<2407.9$ \\
$2408$ & $55.8$ & $212$ & $286.3$ & $572.5$ & $k<596.8$ \\
\bottomrule
\end{tabular}
\end{center}

\bigskip
This completes the proof of Theorem \ref{thm1}.

\textbf{Data availability:} The data mentioned in the paper (computer codes) are available from the fourth author.\\

\paragraph{\textbf{Use of AI Tools:}} An AI assistant (Claude, Anthropic) was used for copyediting and for independent computational verification of the calculations in this paper. The authors take full responsibility for the content.


\end{document}